\documentclass[twoside]{amsart}
\usepackage{amssymb}
\usepackage{amsmath}
\usepackage{amsthm}
\usepackage{graphicx}
\usepackage{color}
\usepackage[all]{xy}

\newtheorem{theorem}{Theorem}[section]
\newtheorem{proposition}[theorem]{Proposition}
\newtheorem{corollary}[theorem]{Corollary}

\newtheorem{lemma}[theorem]{Lemma}

\theoremstyle{definition}
\newtheorem{definition}[theorem]{Definition}

\theoremstyle{remark}
\newtheorem{remark}[theorem]{Remark}

\newcommand{\cO}{{\mathcal O}}
\newcommand{\cI}{{\mathcal J}}

\newcommand{\Q}{\mathbb Q}

\newcommand{\C}{\mathbb C}
\newcommand{\R}{\mathbb R}
\newcommand{\proj}{\mathbb P}

\newcommand{\f}{\varphi}
\newcommand{\ra}{\rightarrow}

\DeclareMathOperator{\codim}{codim}

\begin{document}


\title[Adjunction]{Minimal Model Program with scaling and Adjunction Theory} 
\author[Andreatta]{Marco Andreatta}

\thanks{}

\address{Dipartimento di Matematica, Universit\'a di Trento, I-38050
Povo (TN)} \email{marco.andreatta@unitn.it}

\thanks{Many thanks to Paolo Cascini for explaining some details of the MMP.\\
The author was supported by a grant of the Italian Minister of Research - PRIN} 

\begin{abstract}
Let $(X,L)$ be a quasi-polarized pair, i.e. $X$ is a normal complex projective variety and $L$ is a nef and big line bundle on it. We study, up to birational equivalence, the positivity (nefness) of the adjoint bundles $K_X + rL$ for high rational numbers $r$. For this we run a Minimal Model Program with scaling relative to the divisor $K_X +rL$.
We give then some applications, namely the classification up to birational equivalence of quasi-polarized pairs with sectional genus $0,1$ and of embedded projective varieties $X \subset \proj ^N$ with degree smaller than $2\codim_{\proj ^N} (X) +2$.

\end{abstract}

\maketitle

\section{Introduction}

Let $X$ be a complex projective normal variety of dimension $n$ and $L$ be a nef and big line bundle on $X$. The pair $(X,L)$ is called a {\bf quasi-polarized pair}. The goal of Adjunction Theory is to classify quasi-polarized pairs via the study of the positivity of the adjunction divisors $K_X + rL$, with $r$ positive rational numbers. This has been done extensively in the case in which $L$ is ample, i.e. $(X,L)$ is a polarized pair; \cite{BeltramettiSommese} is the best account on this case. 

However the set up of quasi-polarized pairs is certainly more natural:
in particular when passing to a resolution of the singularities and  taking the pull-back of $L$. The classification of quasi-polarized pairs  will be up to birational equivalence.
Quasi-polarized pairs were first considered by T. Fujita (see \cite{Fujita}). In that paper he made a connection between this theory and the Minimal Model Program (MMP for short); he proved some results under the assumption of the existence of the MMP (more precisely, under the assumption of existence and termination of flips).

The idea of running a MMP supported by a polarizing divisor has been exploited first by M.Reid in the surface case, in the pioneering paper \cite{Re}; his ideas were taken into the case of uniruled $3$-folds by M. Mella in \cite{Me2}.

In this paper, following T. Fujita's ideas as re-proposed by A. H\"oring in \cite{Horing}, and with the use of the MMP developed in \cite{BCHM}, we describe a MMP with scaling related to divisors of type $K_X + rL$ (see Section \ref{MMP}). 

Using the $K_X+rL$-MMP we prove that either the pair $(X,L)$ is birationally equivalent to some very special quasi-polarized
pairs, or it is birationally equivalent to a pair $(X',L')$, which we call a zero-reduction of $(X,L)$,  where $K_{X'}+rL'$ is nef  for $r \geq (n-1)$  (Theorem \ref{teo1}). 

In a further step  we prove that there exists a quasi-polarized pair $(X'',L'')$, which we call a first-reduction of $(X,L)$ and which is related to the original $(X,L)$ via birational equivalences or divisorial contractions to smooth points, such that either $(X'',L'')$ is  in a finite list of special pairs, or $K_{X''} +rL''$ is nef  for $r \geq (n-2)$ (Theorem \ref{teo2}).

Finally we give some applications, namely the classification, up to birational equivalence, of quasi-polarized pairs with sectional genus $0$ or $1$ (Corollary \ref{genus}) and, up to first-reduction,  of embedded projective varieties $X \subset \proj ^N$ with degree smaller than $2\codim_{\proj ^N} (X) +2$  (Corollary \ref{cod}).

\section{Notation and Preliminaries}
Our notation is consistent with the books \cite{BeltramettiSommese} and \cite{KollarMori} and the paper \cite{BCHM}, to which we constantly refer. We give however some basic definitions, for the reader convenience and to state our main objects.

In general $X$ will be a normal, complex and projective variety, that is an irreducible and reduced projective scheme over $\C$, of dimension $n$.
We say that two $\Q$-divisors $D_1, D_2$  are $\Q$-linearly equivalent, $D_1 \sim _{\Q} D_2$, if there exists an integer $m >0$ such that $mD_1$ and  $mD_2$ are linearly equivalent. 
We say that a $\Q$-divisor $D$ is $\Q$-Cartier if some integral multiple is Cartier.  We say that $X$ is $\Q$-factorial if every $\Q$-divisor is  $\Q$-Cartier.

Let $D$ be an  $\R$-divisor; we say that it is nef if $D^.C \geq 0$ for any curve $C \subset X$. We say that it is big if  $D \sim _{\R} A + B$ where $A$ is ample and $B \geq 0$. It is pseudo-effective if it is in the closure of the cone of effective divisors.
Effective or nef divisors are pseudo-effective.

\smallskip
A {\bf quasi-polarized variety} is a pair $(X,L)$ where is $X$ is a complex and projective variety and $L$ is a nef and big Cartier divisor (equivalently a nef and big line bundle). 

\smallskip
A log pair $(X, \Delta)$ is a normal variety  $X$ and an effective $\R$ divisor $\Delta$ such that $K_X + \Delta$ is $\R$-Cartier. 
 A log resolution of the pair $(X,\Delta)$ is a projective birational morphism $g: Y \ra X$ such that $Y$ is smooth and the exceptional locus is a divisor which, together with $g^{-1}(\Delta)$, is simple normal crossing. We can write
 
 $$g^*(K_X + \Delta) = K_Y + \Sigma b_i \Gamma_i , $$ 
 
 where $\Gamma _i$ are distinct prime divisors.
 
 The log pair $(X,\Delta)$ is Kawamata log terminal (klt) if for every (equivalently for one) log resolution as above $b_i  < 1$, for all $i$. If $\Delta=0$ and $b_i <0$, for all $i$,  then $X$ has terminal singularities.

\section{Quasi-polarized Pairs and Adjunction Theory}

\begin{definition} \label{bireq}
Two quasi-polarized pairs $(X_1,L_1)$ and $(X_2,L_2)$ are said to be {\bf birationally equivalent} if there is another variety $Y$ with birational morphisms $\f_i :Y \ra X_i$ such that $\f_1^*L_1 = \f_2^*L_2$.
\end{definition}

\smallskip
\begin{definition} \label{defgenus}
The Hilbert polynomial of the quasi-polarized pair $(X,L)$ is given by $\chi(X,tL) = \Sigma  _{j=0,....,n} \chi_j t^{[j]} /j !$, where $t^{[j]} =  t(t+1)...(t+j-1) $, $t^{[0]} = 1$ and  $\chi_0, ..., \chi_n$ are integers. 

By the Riemann-Roch Theorem we have that $\chi_n = L^n$ and, if $X$ is normal, that $-2 \chi_{n-1} = (K_X +(n-1) L)^.L^{n-1}$, for a canonical divisor $K_X$ on $X$.

The sectional genus of the pair $(X,L)$ is defined as $g(X,L) = 1- \chi_{n-1}$. 

The $\Delta$-genus, on the other hand,  is defined as $\Delta(X,L) = n + \chi_n - h^0(X,L)$.
\end{definition}
 
\smallskip
Assume that $X$ has at most terminal singularities and that  $K_X$ is not nef. In this paper, in the spirit of Mori theory, an {\bf extremal ray} $R$ on $X$ will be an extremal ray in the cone $\overline{NE(X)} \subset N^1(X)$ which has negative intersection with the canonical class. In particular, by a theorem of Mori,  $R= \R^+[C]$, where $C\subset X$ is a rational curve such that $-{K_X}^.C >0$. 

Let $\f_R: X \ra Z$ be the contraction associated with the extremal ray $R$: that is $\f_R$ is a morphism with connected fibers onto a normal projective variety $Z$ and a curve in $X$ is in a fiber of $\f_R$ if and only if its class is in the ray $R$. 
The existence of $\f_R$ is the famous base-point-free theorem of Kawamata-Shokurov in the MMP theory.

If $\f_R$ is of fiber type, i.e. $\dim X > \dim Z$, then  $\f_R: X \ra Z$ is called a {\bf Mori fiber space}.

Otherwise the contraction
$\f_R$ is birational and it can be either divisorial or small.

\smallskip
For a normal quasi-polarized pair $(X,L)$ let 
$$r(X,L):=\hbox{sup} \{t \in \R : tK_X + L \hbox{ is nef}\}.$$
By the rationality theorem of Kawamata, $r(X,L)$ is a rational non negative number. 
If $r(X,L) \not= 0$ we define $\tau(X,L) := 1/ r(X,L)$.

\begin{lemma} \label{amplefiber} 
Let $\f_R: X\ra Z$ be a Mori fiber space associated with the extremal ray $R= \R^+[C]$.
A nef and big  line bundle $L$ on $X$ is 
$\f_R$-ample, i.e. $L^.C >0$.  
\end{lemma}

\proof If, by contradiction, $L^.C=0$ then there exists a line bundle $A$ on $Z$ such that $L= \f_R^* (A)$
(see Corollary 3.17 in \cite{KollarMori}).
This implies that $L^n =0$, which is a contradiction since $L$ is nef and big (see Proposition 2.61 in \cite{KollarMori}).

\begin{remark} \label{ample}
Let $\f_R: X\ra Z$ be the contraction associated with the extremal ray $R= \R^+[C]$.
Assume that $L$ is a $\f$-ample line bundle. By adding to $L$  the pull-back of a sufficiently ample line  bundle from $Z$
we can assume that 

i) $L$ is ample, 

ii)  $r(X,L) \not= 0$ and

iii) the intersections of $(K_X + \tau(X,L) L)$ with curves whose classes are in $R$ are zero and they are positive with all other curves.

(The proof of these remarks is standard in the theory of ample line bundles; use for instance section 1.5, in particular Proposition 1.45,  of \cite{KollarMori}).
\end{remark}

\begin{proposition} \label{Morifiberspace} Let $X$ be a normal variety with terminal $\Q$-factorial singularities and
$L$ be a nef and big line bundle on $X$.
Let $\f_R: X\ra Z$ be a Mori fiber space associated with the extremal ray $R= \R^+[C]$ and 
$r$ be a positive rational number such that  $(K_X+ rL)^.C < 0$. Note that this implies that 
$\tau(X,L) > r$ (possibly adding to $L$ the pull-back of a sufficiently ample line  bundle from $Z$).

A) If $r \geq (n-1)$ then  $(X,L)$ is one of the following pairs: 
\begin{itemize}
\item $(\proj ^n, \cO(1))$ and $r < (n+1)$,
\item $(Q, \cO(1)_{|Q})$, where $Q\subset \proj ^{n+1}$ is a quadric and $r < n$,
\item $C_n(\proj^2, \cO(2))$, a generalized cone over $(\proj^2, \cO(2))$, and $r < n$,
\item $\f_R$ gives to $X$ the structure of a  $\proj^{n-1}$-bundle over a smooth curve $C$ and $L$ restricted to any fiber is $\cO(1)$ and $r < n$.
\end{itemize}

B) If $r \geq (n-2)$ then  $(X,L)$ is one of the following pairs:
\begin{itemize}
\item one of the pair in the previous list,   
\item a del Pezzo variety, that is $-K_X  \sim _{\Q} (n-1)L$ with $L$ ample,  $r < (n-1)$, 
\item $(\proj ^4, \cO(2))$,
\item $(\proj ^3, \cO(3))$,
\item $(Q, \cO(2)_{|Q})$, where $Q\subset \proj ^4$ is a quadric,
\item $\f_R$ gives to $X$ the structure of a  quadric fibration over a smooth curve $C$ and $L$ restricted to any fiber is 
$\cO(1)_{|Q}$, $r < (n-1)$,
\item $\f_R$ gives to $X$ the structure of a $\proj^{n-2}$-bundle over a normal surface $S$ and $L$ restricted to any fiber is $\cO(1)$, $r < (n-1)$,
\item $n=3$, $Z$ is a smooth curve, the general fiber of $\f_R$ is $\proj ^2$ and $L$ restricted to it is $\cO(2)$.
\end{itemize}

\end{proposition}

\proof Using Lemma \ref{amplefiber} and Remark \ref{ample} we can assume that $L$ is ample. The Proposition follows then by the 
''classic'' adjunction theory developed by T. Fujita and by A. Sommese and his school: more precisely the results are summarized in section 7.2 and 7.3 of \cite{BeltramettiSommese}. One of my personal contribution to this theory is its extension to the case with terminal or even log-terminal singularities in the papers \cite{An1} and \cite{An2}

\begin{proposition} \label{birational}  Let $X$ be a normal variety with terminal $\Q$-factorial singularities and
$L$ be a nef and big line bundle on $X$.
Let $\f_R: X\ra Z$ be a birational contraction associated with the extremal ray $R= \R^+[C]$ and let $r$ be a rational number such that  $(K_X+ rL)^.C  < 0$.

Assume that $r \geq (n-2)$ and that  $L^.C \not= 0$.

Then $r < (n-1)$ and $\f_R$ is the contraction of an irreducible divisor $E$ into a smooth point $z\in  Z$. 

Therefore $L':= \f_*(L)$ is a nef and big Cartier divisor (i.e. line bundle) on $Z$ and  $\f_R^*L' = L + bE$ with $b > 0$.
\end{proposition}

\proof By the assumption, using Remark \ref{ample}, we can assume that $L$ is ample, that  the intersection of $(K_X + \tau(X,L) L)$ with curves whose classes are in $R$ is zero and positive on all other curves and that $\tau:= \tau(X,L) > r \geq (n-2)$.

We can apply Theorem 2.1 of \cite{An2} and have that, if $F$ is a component of a non trivial fiber of $\f_R$, then $\dim F \geq \tau > r$. Since $r\geq (n-2)$ it implies that $F$ is a divisor and $r < n-1$. Therefore $\f_R$ is divisorial with exceptional locus an irreducible divisor $E = F$;  in particular $Z$ has terminal $\Q$-factorial singularities
(see Proposition 3.36 and Corollary 3.43 in \cite{KollarMori}). 

Let us prove that $\f_R(E) :=z$ is smooth; we will do it by induction on $n$, the dimension of $X$.  Let $n=2$; terminal singularities in dimension $2$ are smooth and  a birational extremal contraction on a smooth surface  is the contraction of a $-1$-curve.
By Castelnuovo theorem the image is a smooth point.

Let $n >2$. The problem is local and therefore we can assume that $Z$ is affine,
i.e. $\f_R: X \ra Z$ is a Fano-Mori contractions supported by $K_X + \tau L$ as in the set up described in section 2 of \cite{AW1}. By the main Theorem of \cite{AW1} we can assume that 
$L$ is spanned by global section (or even very ample: see Propositions 1.3.3 and 1.3.4 in \cite{AW2}). Take a general section in the linear system $X' \in |L|$; $X'$ has dimension $(n-1)$ and, by Bertini theorem, it has terminal $\Q$-factorial singularities. Moreover $(f_R)_{|X'}: X' \ra Z'$ is a divisorial Fano Mori contraction supported by  $K_{X'} + (\tau -1) L_{|X'}$
(see Proposition 2.6  in \cite{AW1}).

By the inductive assumption $z\in Z'$ is smooth. 
The smoothness of $z\in Z$  follows then from the next lemma, which is Lemma 1.7 in \cite{Me} (the case $n=3$ is in the proof of Lemma 5.3 in \cite{Kollaral}).

\begin{lemma} Let $Z$ be a variety with terminal $\Q$-factorial singularities and $Z'$ a smooth divisor on $Z$. Then $Z'$ is Cartier; in particular, $Z$ is smooth along $Z'$.
\end{lemma}

Since $L$ is a Cartier divisor,  $L':= \f_*(L)$ is a Cartier divisor outside $z$; since $z$ is smooth $L'$ is Cartier everywhere. 

Let $\f_R^*L' = L + bE$; since $L^.C >0$ we have that $b>0$.

\begin{corollary} \label{Kawakita}
If $n=3$ the divisorial contractions of Proposition \ref{birational} are weighted blow-ups of smooth points with weights $(1,1,b)$, with $b \geq 1$ an integer number.
\end{corollary}

\proof Divisorial contractions on terminal $\Q$-factorial $3$-folds which
send a divisor into smooth points are classified in \cite{Kaw}: they are weighted blow-ups with weights $(1,a,b)$, with $(a,b) =1$.   Let $l$ be the curve on $E = \proj(1,a,b)$ numerically equivalent to the class $\cO(1)$. We have that $- {K_X} ^.l = (a+b)/ab$ and, under our assumption,  also that $- {K_X }^.l = \tau(X,L) L^.l > 1$. This implies that $a=1$ and the nef value is $1 <(b+1)/b \leq 2$. 

\medskip
\begin{remark}\label{factorial} 
In a forthcoming paper, together with L. Tasin, we are going to extend
the above Corollary \ref{Kawakita} for any $n\geq 3$: i.e. the divisorial contractions of Proposition \ref{birational} are weighted blow-ups of smooth points with weights $(1,1,b, ..., b)$, with $b \geq 1$ an integer number.
Note that if $X$ is factorial then $b=1$, i.e. in the factorial case these maps are the simple blow-ups of smooth points. 
\end{remark}


\section{Minimal Model Program with scaling} 
\label{MMP}

Let $(X,L)$ be a quasi-polarized variety and assume that $X$ has at most terminal $\Q$-factorial singularities.
Let also $r$ be a positive rational number. 

\begin{lemma} \label{asymptotic}
Under the above assumption (in particular  with $L$ nef and big) there exists an effective $\Q$-divisor $\Delta^r$ on $X$ such that

$$\Delta^r \sim_{\Q} r L \hbox{\ \ \ and\ \ \  }
(X, \Delta^r) \hbox{\ \ is Kawamata log terminal.}$$
\end{lemma}

\proof This lemma is well known to specialists and it can be proved in different ways. Since $L$ is nef and big 
the asymptotic multiplier ideal of $rL$ is trivial, i.e. $\cI(X,||rL||)= \cO_X$ (Proposition 11.2.18 in \cite{Laz} in the smooth case
or Corollary 5.2 in \cite{CaBi} in the terminal case and under the weaker assumption that $L$ is nef and abundant).  Take $D$ a generic divisor in $mrL$ for sufficiently large $m$ and let $\Delta^r := {1\over m}D$.
$\Delta^r$ is effective and $\Q$ linearly equivalent to $rL$.

 Moreover, for $m$ sufficiently large,
$\cI(X,\Delta^r)= \cI(X,{1\over m}(|mrL|) = \cI(X,||rL||)= \cO_X$, i.e. $(X, \Delta^r)$  is Kawamata log terminal.

\medskip
Consider the pair $(X, \Delta^r)$ and the $\Q$-Cartier divisor $K_X + \Delta^r  \sim_{\Q} K_X + r L$.

By Theorem 1.2 and Corollary 1.3.3 of \cite{BCHM}  we can run a 

\centerline{\bf{$K_X + \Delta^r$- Minimal Model Program with scaling:}}

\smallskip
\centerline{$(X_0, \Delta_0^r) = (X, \Delta^r) \ra (X_1, \Delta_1^r)   \ra ----\ra (X_s, \Delta_s^r)$}

such that:

1) the log pair $(X_i,  \Delta_i^r)$ is a Kawamata log terminal, for $i =0, ..., s$. 

2) Each map $\f_i : X_i \ra X_{i+1}$ is a birational map which is either a divisorial contraction 
or a flip associated with an extremal ray $R_i$; in particular $X_i$ has at most terminal $\Q$-factorial singularities, for $i =0, ..., s$.

3) a) If $K_X + \Delta^r$ is pseudo-effective then $K_{X_s} + \Delta_s^r$ is nef, 

\ \ \ \  b)  if  $K_X + \Delta^r$ is not pseudo-effective then  $X_s$ is a Mori fiber space relatively to $K_{X_s} + \Delta_s^r$.

\medskip
The next Proposition has been proved  by T. Fujita in section 4 of \cite{Fujita},  under the assumption of the existence of  Minimal Models (more precisely subordinated to the Flip conjecture). A.  H\"oring (\cite{Horing}) has adapted Fujita's argument to the notations and the spirit of \cite{BCHM}; see the Claim in the course of the proof of his Proposition 1.3.

\begin{proposition}  \label{rL} Under the above notation and assumptions
suppose moreover that $r \geq (n-1)$. 

For every $i = 0, ..., s$, we have ${\Delta_i^r} ^. R_i = 0$; therefore there exist nef and big Cartier 
divisors $L_i$ on $X_i$ such that $\f_i^*(L_{i+1}) = L_i$ and $\Delta_i^r \sim_{\Q} r L_i$.  

Thus at every step of the above given MMP we have a quasi-polarized variety $(X_i, L_i)$, with at most terminal $\Q$-factorial singularities.

Note also that $H^0(K_{X_i} +tL_i) = H^0(K_{X_{i+1}} +tL_{i+1})$ for any $t\geq 0$.
\end{proposition}

\proof 
The Proposition follows by induction on $i$. 

Each map $\f_i : X_i \ra X_{i+1}$ is a birational map associated with an extremal ray $R_i=\R^+[C_i]$ with 
$$(K_{X_i} + r L_i)^.C_i = (K_{X_i} + \Delta_i^r)^.C_i  < 0.$$ 

Since $r \geq (n-1)$, by Proposition \ref{birational}, we have that ${L_i} ^. C_i = {\Delta_i^r} ^. C_i = 0$. 

Let $\f_{R_i}: X_i \ra Z$ be the contraction of the extremal ray $R_i$; since ${L_i} ^. C_i = 0$ there exists a 
nef and big line bundle $L$ on $Z$ such that $\f_{R_i}^* L = L_i$ (Corollary 3.17 of \cite{KollarMori}). 

If $\f_{R_i}$ is birational then $\f_{R_i}= \f_i$ and we take $L_{i+1}$ to be $L$ itself.  If $\f_{R_i}$ is small let $\f^+: X_{i+1} \ra Z$ be its flip; define then  $L_{i+1}$  to be ${\f^+} ^*(L)$.

Note that $\Delta_{i+1}^r = \f_*\Delta_i^r  \sim_{\Q}  \f_*(rL_i)= rL_{i+1}$.

Let us prove the last statement, namely $H^0(K_{X_i} +tL_i) = H^0(K_{X_{i+1}} +tL_{i+1})$  for any $t\geq  0$.
This is obvious if $\f_i$ is a flip, since $X_i$ and $X_{i+1}$ are isomorphic in codimension $1$. If $\f_i$ is birational
then $\f_i^*(K_{X_{i+1}}) + \alpha_i E_i = K_{X_i}$, where $E_i$ is the effective exceptional divisor and $\alpha_i >0$;
in particular $\f_i^*(K_{X_{i+1}} +tL_{i+1})  + \alpha_i E_i = K_{X_i} +tL_i$. Our claim follows applying for instance Lemma 7.11 in \cite{deb}.

\begin{remark} 
\label{n-1}

Proposition \ref{rL} says that, in a 
$K_X + \Delta^r$- Minimal Model Program  with scaling with $r \geq (n-1)$, the birational contractions $\f_i$ are trivial with respect to the boundaries  $\Delta_i^r \sim_{\Q} r L_i$.
This implies that the birational part of these $K_X + \Delta^r$- Minimal Model Programs is the same for any $r\geq (n-1)$.

If the last pair of the program, $(X_s, L_s)$, has a Mori fiber space structure relative to the ray $R= \R^+[C]$ , then $\Delta_s^r$ is not trivial on it, i.e. $(\Delta_s^r)^.C >0$ (see \ref{amplefiber}). 

Therefore  $(K_{X_s}+ {\Delta_s^{r}})^. C < 0$ but,  for $r' > r$, the divisor $K_{X_s}+ \Delta_s^{r'}$ can be nef.

\end{remark}

\begin{definition}  \label{def1}Let $(X,L)$ be a quasi-polarized variety such that $X$ has at most terminal $\Q$-factorial singularities.
Let $(X_s,L_s)$ be a quasi-polarized pair such that $X_s$ is the last variety in a $K_X + \Delta^{(n-1)}$- Minimal Model Program  with scaling and $L_s$ is the corresponding nef and big line bundle on $X_s$ coming from Proposition \ref{rL}.  

 We will  call  $(X',L') := (X_s,L_s)$ a    {\bf zero-reduction} of the pair $(X,L)$.
\end{definition}

\begin{remark}  
i) A zero-reduction is birationally equivalent to the original pair.

ii) Long ago with A. Sommese we studied the surface case ($n=2)$  in \cite{AnSo}; in particular Proposition 1.7 in that paper gives the construction of the zero-reduction for Gorenstein surfaces (note that for $n=2$ terminal singularities are actually smooth).

\end{remark}


\section{Adjunction theory via MMP with scaling}

The results of this section are classical in the case when $L$ ample
and $X$ is smooth; \cite{BeltramettiSommese} is the best reference for that. 
They have been proved recently in \cite{BKLN}, with $X$ smooth and $L$ nef, but with 
the extra assumption that $L$ is strictly positive (ample) on the negative part of the Mori cone,
i.e. $\{ C \hbox{ curve in } X: L^.C=0\} \subset \{ C \hbox{ curve in } X: {K_X}^.C \geq 0\} $.

Beside the large generality in which we prove them, i.e. $X$ has terminal $\Q$-factorial singularities and $L$  is
nef and big, we would like to remark that our approach is different: i.e. we use the existence of the MMP
and then we apply adjunction techniques to describe the steps of the program  in details. 

Most of the results would actually work also in the case in which $X$ has log terminal $\Q$-factorial singularities.

\subsection{Adjunction on the zero-reduction} The first step in the Adjunction Theory of quasi-polarized pairs
is given by the following theorem; 
Part 3)  was first proved by A. H\"oring  (\cite{Horing}, Proposition 1.3).

\begin{theorem} \label{teo1}
Let $(X,L)$ be a quasi-polarized variety  such that  $X$ has at most terminal $\Q$-factorial singularities.

1) $K_{X} + (n+1)L$  is pseudo-effective and on a zero-reduction $(X',L')$  the $\Q$-Cartier divisor $K_{X'} + (n+1)L'$ is nef.

2) $K_{X} + nL$  is not pseudo-effective if and only if any zero-reduction $(X',L')$ is $(\proj ^n, \cO(1))$. 

If  $K_{X} + nL$  is pseudo-effective then on a zero-reduction $(X',L')$  the $\Q$-Cartier divisor $K_{X'} + nL'$ is nef.

3) $K_{X} + (n-1)L$ is not pseudo-effective if and only if any zero-reduction $(X',L')$ is one of the pairs in \ref{Morifiberspace} A). 

If $K_{X} + (n-1)L$ is pseudo-effective then on a zero-reduction $(X',L')$  the $\Q$-Cartier divisor $K_{X'} + (n-1)L'$ is nef.

\end{theorem}

\proof We use  the construction in Section \ref{MMP} and the Proposition \ref{Morifiberspace}. 

Run a $K_X + (n+1)L $-Minimal Model Program  on $(X,L)$ as in Section \ref{MMP} and let $(X_s,L_s)$ be the last pair of the process;  by Corollary \ref{n-1},  $(X',L') := (X_s,L_s)$ is
a zero-reduction of the pair $(X,L)$, as in Definition \ref{def1}).   

Assume, by contradiction, that $K_{X} + (n+1)L$  is not pseudo-effective; then $X_s$ is a Mori fiber space associated with an extremal ray $R= \R^+[C]$
such that $(K_{X_s} + (n+1) L_s ) ^. C < 0$. This cannot be the case, by Proposition \ref{Morifiberspace},
which in fact says that if $(K_{X_s} + r L_s ) ^. C < 0$ then $r < (n+1)$. 

Therefore  $K_{X} + (n+1)L$  has to be pseudo-effective and, on  a zero-reduction $(X',L')$, the divisor
$K_{X'} + (n+1) L'  \sim_{\Q} K_{X_s}+ \Delta_s^{n+1}$ is nef.

\smallskip
Points 2) and 3) can be proved similarly; let us prove for instance point 3). 
Let $(X',L')$ be a zero-reduction of $(X,L)$ defined in \ref{def1}: if $K_{X} + (n-1)L$ is not pseudo-effective
then $X'$ is a
Mori fiber space associated with an extremal ray $R= \R^+[C]$
such that $(K_{X'} + (n-1) L' ) ^. C < 0$.  By Proposition \ref{Morifiberspace}, it has to be one of the pairs in \ref{Morifiberspace} A). 

If  $K_{X} + (n-1)L$ is pseudo-effective  then $K_{X'} + (n-1)L' $ is nef.

\begin{corollary}
\label{alpha}
On the zero-reduction $(X',L')$   there are no extremal rays $R = \R^+[C]$ such that ${L'}^.C =0$.
\end{corollary}

\proof In fact $K_{X'} + (n+1)L'$ is nef and therefore for every curve $C \subset X'$ such that $-{K_{X'}} ^. C < 0$ it must be 
${L'} ^. C >0$.

\begin{remark} The zero-reduction is related to the almost holomorphic map constructed in \cite{CaPe}, a reduction map
for  nef line bundles. Actually their map will factor through the zero-reduction; this last in fact 'reduces'  the curves on which $L$ is zero and whose classes are in extremal rays. 
In particular on a Fano variety $X$ a zero-reduction is a map as in \cite{CaPe}.
\end{remark}


\subsection{First-reduction for quasi-polarized pairs}

 Let $(X,L)$ be a quasi-polarized pair with at most terminal $\Q$-factorial singularities and 
let $(X',L')$ be the zero-reduction of $(X,L)$.  
We proceed now with a further step in Adjunction theory. 

\medskip
Let $r \geq (n-2)$ and, as in Lemma \ref{asymptotic}, 
take an effective $\Q$-divisor ${\Delta'}^r$  on $X$ such that:
${\Delta'}^r \sim_{\Q} r L'$ and $(X', {\Delta'}^r)$ is Kawamata log terminal.

Consider a 
$K_{X'} + {\Delta'}^r$- Minimal Model Program  with scaling as in the first part of Section \ref{MMP}. 

\centerline{$(X'_0, {\Delta'}_0^r) = (X', {\Delta'}^r) \ra (X'_1, {\Delta'}_1^r)   \ra ----\ra (X'_s, {\Delta'}_s^r)$}

\begin{proposition} \label{frL}
Under the above notations and assumptions, at every step $i = 0, ..., s$, the morphism $\f'_i: X'_i\ra X'_{i+1}$  is the contraction of an irreducible divisor to a smooth point; in particular
$X'_{i+1}$ has at most terminal $\Q$-factorial singularities.  

On $X'_{i+1}$ there exists a nef and big line bundle $L'_{i+1}$ such that 
 ${\f'}_i^*(L'_{i+1}) = L'_i+ bE_i$, where $E_i$ is the exceptional divisor and 
 $b$ a positive integer.
 
 Thus at every step of the above given MMP we have a quasi-polarized variety $(X'_i, L'_i)$, with at most terminal $\Q$-factorial singularities; moreover  ${\Delta'}_i^r \sim_{\Q} r L'_i$.

We also have $H^0(K_{X'_i} +tL'_i) = H^0(K'_{X_{i+1}} +tL'_{i+1})$, for any $t= 0, ...., r$.

\end{proposition}

\proof 
The proof is by induction on $i$. 
Each map $\f'_i : X_i \ra X_{i+1}$ is a birational map associated with an extremal ray $R_i=\R^+[C_i]$ such that 
$$(K_{X_i} + r L_i)^.C_i = (K_{X_i} + \Delta_i)^.C_i  < 0.$$ The Proposition will follow directly from 
Proposition \ref{birational} if we prove that ${\Delta'_i} ^. R_i = r{L'_i}^. R_i \not=0$.

By Corollary \ref{alpha} this is the case for $i= 0$.
Assume by contradiction that, at a further step $k$, we encounter, for the first time, a ray $R_k= \R^+ [C_k]$ with ${L'_k}^. C_k = 0$. 
At the previous step, by induction, $\f'_{k-1} : X'_{k-1} \ra X'_k$ is the contraction of an irreducible divisor to a smooth point $p$, associated with the ray $R_{k-1}$,  and $L'_{k-1} = {\f'}_{k-1} ^* L'_k + bE_k$, with $b > 0$. Therefore we have that ${L'_{k-1}} ^. \bar C = - b{E_k} ^. \bar C$, where $\bar C$ is the strict transform of $C_k$. Since $L'_{k-1}$ is nef and $E_k$ effective this implies that these intersections are zero and $\bar C$ doesn't pass through $p$. 
The composition $\f'_k \circ \f'_{k-1}$ is a contraction of an extremal face of dimension two, i.e. it is generated by two extremal rays; one is $R_{k-1}$. The ray $\R^+[\bar C] $ is contracted by $\f'_k \circ \f'_{k-1}$ and its locus  is disjoint from the one of $R_{k-1}$. Therefore $\R^+[\bar C] $ is the other extremal ray of this face on $X'_{k-1}$. But it has zero intersection with $L'_{k-1}$, which is a contradiction.

As for the claim that ${\Delta'}_i^r \sim_{\Q} r L'_i$, recall that  ${\Delta'}_{i+1}^r := \f'_*{\Delta'}_i^r$ and  $rL'_{i+1}=  (\f'_*(rL'_i))^{**}$.
By the inductive assumption, we can assume that ${\Delta'}_i^r \sim_{\Q}  rL'_i$, i.e. that there exists an integer $m$ such that $m{\Delta'}_i^r \in |mrL'_i|$. Since $\f'_i$ is the contraction of a divisor to a smooth point, $m{\Delta'}_{i+1}^r$ is Cartier and  $m{\Delta'}_{i+1}^r \in |mrL'_{i+1}|$, i.e. 
 ${\Delta'}_{i+1}^r \sim_{\Q}  rL'_{i+1}$.
 
\smallskip
The proof that $H^0(K_{X_i} +tL_i) = H^0(K_{X_{i+1}} +tL_{i+1})$, for any $t= 0, ...., r$, is similar to the one in 
Proposition \ref{rL}, using for instance Lemma 7.11 in \cite{deb}.

\begin{definition} Let $(X,L)$ be a quasi-polarized variety such that $X$ has at most terminal $\Q$-factorial singularities
and let $(X',L')$ be  a zero-reduction.
Let $(X'_s,L'_s)$ be a quasi-polarized pair where $X_s$ is the last variety in a $K_X' + \Delta'^{(n-2)}$- Minimal Model Program and $L'_s$ is the corresponding nef and big line bundle on $X'_s$ coming from Proposition \ref{frL}.  Let 
$\rho: X' \ra X''$ be the composition $\rho = \f'_{s-1} \circ ... \circ\f'_0$

 We will  call  $(X'',L'') = (X'_s,L'_s)$, together with a zero-reduction $X'$ and the map $\rho: X' \ra X''$, a  {\bf first-reduction of the pair $(X,L)$}.
\end{definition}

\begin{remark}
i) In the case $L$ is ample and $X$ is factorial the definition of first-reduction is in agreement with the Sommese's definition, i.e. Definition 7.3.3 in \cite{BeltramettiSommese} (see Remark \ref{factorial}).

In the case $n=3$, with the additional assumption that $L$ is movable and without base points, the first-reduction was studied in \cite{Me2} (Corollary 3.10 in that paper).

\smallskip
ii) The pairs $(X,L)$ and $(X'',L'')$ are not birationally equivalent. 

However the morphism $\rho: X' \ra X''$ is very simple, namely it consists of a series of divisorial contractions to smooth points. By Corollary \ref{Kawakita} (if $n=3$, in general by Remark \ref{factorial}), they are weighted blow-ups of weights $(1,1,b, ..., b)$, with $b \geq 1$. 

\smallskip
iii) We could run directly a $K_X + \Delta^{(n-2)}$- Minimal Model Program  with scaling on $(X,L)$.
In this case, with the help of Proposition \ref{birational}, we have  at each step $i$ two possibilities, namely either ${\Delta_i} ^. R_i = 0$ or ${\Delta_i} ^. R_i \not= 0$. In the first case we define a nef and big line bundle on $L_{i+1}$ on $X_{i+1}$
such that  $\f_i^*(L_{i+1}) = L_i$. In the second case, $\f_i: X_i\ra X_{i+1}$  is the contraction of an irreducible divisor to a smooth point; then we can define a
nef and big line bundle $L_{i+1}$ on $X_{i+1}$,  such that 
 $\f_i^*(L_{i+1}) = L_i+ b_jE_i$ ($E_i$ is the exceptional divisor of the blow-up).
 
At the end we will reach a quasi-polarized pair $(X_s, L_s)$ which has the same property of the first-reduction $(X'',L'')$
in Theorem \ref{teo2}.

The above construction, which splits the Program in two parts, namely a first part contracting all rays whose intersection with the polarization is zero and a second consisting of weighted blow-ups of smooth points, is more accurate and useful.
 \end{remark}

\smallskip
Using the first-reduction we can push adjunction theory a step further.

\begin{theorem} \label{teo2}
Let $(X,L)$  be a quasi-polarized variety  such that  $X$ has at most terminal $\Q$-factorial singularities.

1) $K_{X} + (n-2)L$ is not pseudo-effective if and only if any first-reduction  $(X'', L'')$ is one of the pairs in  \ref{Morifiberspace} A) or B). 

2) If $K_{X} + (n-2)L$ is pseudo-effective then on any first-reduction $(X'', L'')$  the divisor 
$K_{X''} + (n-2)L''$ is nef.
\end{theorem}

\proof The proof is similar to the one of Theorem \ref{teo1}. Take  a $K_{X'} + {\Delta'}^{(n-2)}$- Minimal Model Program ending in the first-reduction $(X'', L'')$. If $K_X +(n-2)L$ is not pseudo-effective  then $(X'', L'')$ is a Mori fiber space and, by Proposition \ref{Morifiberspace},  we are as in point 1).  Otherwise $K_{X''} + (n-2)L''$ is nef.


\section{Applications}

\medskip
The next corollary was proved by T. Fujita (\cite{Fujita}), under the assumption of the existence of a Minimal Model for $X$. After \cite{BCHM} parts 1) and 2) has been first proved in (\cite{Horing}) and part 3) in \cite{Fukuma}.

\begin{corollary} \label{genus}
Let $(X,L)$ be a quasi-polarized variety and $g(X,L)$ be its sectional genus (see Definition \ref{defgenus}). Then:

1)  $g(X,L) \geq 0$.

2)   $g(X,L) = 0$ if and only if $(X,L)$
is birationally equivalent to one of the following quasi-polarized pairs:
\begin{itemize}
\item $(\proj ^n, \cO(1))$, or 
\item $(Q, \cO(1)_{|Q})$, where $Q\subset \proj ^{n+1}$ is a quadric, or 
\item$C_n(\proj^2, \cO(2))$, a generalized cone over $(\proj^2, \cO(2))$,    or
\item $X$ has the structure of a  $\proj^{n-1}$-bundle over a smooth rational curve $C$ and $L$ restricted to any fiber is $\cO(1)$ (a scroll over a rational curve).
\end{itemize}

3) )If  $X$ is normal then $g(X,L) = 1$ if and only if $(X,L)$
is birationally equivalent to one of the following quasi-polarized pairs:
\begin{itemize}  
\item a del Pezzo variety, i.e. $-K_{X'}  \sim _{\Q} (n-1)L'$ with $L'$ ample, 
\item $X'$ has the structure of a  $\proj^{n-1}$-bundle over an elliptlic curve $C$ and $L'$ restricted to any fiber is $\cO(1)$ (a scroll over an elliptic curve).
\end{itemize}
\end{corollary}

\proof Let $\nu: X' \ra X$ be the normalization of $X$; it is straightforward to see that $g(X',\nu^*L) \leq g(X,L) $ (see for instance \cite{Horing}, p. 128). Therefore in proving points 1) and 2) we can assume that $X$ is normal.

By Lemma 1.8 in \cite{Fujita} the sectional genus is a birational invariant of normal quasi-polarized pairs; so we can replace $(X,L)$ first with its resolution and then with its zero-reduction. Call this new pair $(X',L')$.

By Theorem \ref{teo1} if $K_{X'} + (n-1)L'$ is not nef then $(X',L')$ is one of the pair in \ref{Morifiberspace} A). They give the first three cases in 2) and the case in which $(X,L)$ is a scroll over a curve $C$. In this last  case $g(X,L) =g(C)$ and we get the fourth case in 2) and the second in 3).

We can thus assume that $K_{X'} + (n-1)L'$ is nef; therefore  $2g(X', L') -2 =(K_{X'} + (n-1)L' ){L'} ^{n-1} \geq 0$, i.e. $g(X,L) \geq 1$. 

Assume that $g(X', L') =1$. By the previous equality, the facts that  $(K_{X'} + (n-1)L' )$  is nef and $L'$ is nef and big, we get that  $K_{X'} + (n-1)L' $ is numerically trivial. It is straightforward to see that $K_{X'} + (n-1)L' $ is effective
(see for instance \cite{Fujita}, p. 115). Therefore  $K_{X'} + (n-1)L' $ is trivial and we are in the first case of 3).

\medskip
The following application extends the main result in \cite{Ionescu} from the case of smooth embedded varieties to the singular ones. The extension to the special case in which $X$ has at most crepant singularities, has been given in \cite{BKLN}, Theorem 4.2.

\begin{corollary}\label{cod}
 Let $X \subset \proj^N$ be a non degenerate projective variety of dimension $n \geq 3$ and of degree $d$. 
Assume that $d < 2\codim_{\proj ^N} (X) + 2$ (equivalently that $d > 2 \Delta(X,\cO_X(1))$, where $\Delta$ is the
delta genus of the pair, see Definition \ref{genus}). 
Then either $(X,\cO(1))$ is birationally equivalent to one of the quasi-polarized pair in Proposition \ref{Morifiberspace} A)
or the first-reduction of the resolution of $X$ is one of the quasi-polarized varieties in Proposition \ref{Morifiberspace} B).
\end{corollary}

\proof Let $\pi: \tilde X \ra X$ be a resolution of the singularities of $X$: let also $\tilde L:= \pi^* L$. 
$\tilde L$ is globally generated and $h^0(\tilde X, \tilde L) \geq N+1$. The Corollary follows
by Theorem \ref{teo2} 1) if we show that $K_{\tilde X} + (n-2) \tilde L$ is not pseudo-effective. 

Take $L_1, ..., L_{n-1}$ general members in $|\tilde L|$ and let  
$C := L_1 \cap ...\cap L_{n-1}$;
Lemma A.2 in \cite{Me2} says that 
$(K_{\tilde X} + (n-2) \tilde L) ^. C  < 0$.
By Theorem 0.2 in \cite{BDPP} this implies that $K_{\tilde X} + (n-2) \tilde L$ is not pseudo-effective.

\end{document}